\documentclass{article}
\usepackage{amsmath, amssymb}
\newcommand{\sqtimes}{\boxtimes}
\newcommand{\la}[1]{\mathfrak{#1}}

\newtheorem{thm}{Theorem}[section]
\newtheorem{dfn}[thm]{Definition}
\newtheorem{lem}[thm]{Lemma}
\newtheorem{prop}[thm]{Proposition}
\newtheorem{cor}[thm]{Corollary}

\newcommand{\C}{{\mathbb C}}
\newcommand{\R}{{\mathbb R}}
\newcommand{\Z}{{\mathbb Z}}

\title{Towards Vertex Algebras of Krichever-Novikov Type, Part I}
\author{\textsc{Klaus J. Linde}
\footnote{\textsl{Mathematisches Institut, Ludwig-Maximilian Universit\"{a}t, Theresienstr. 39, 80333 M\"{u}nchen, Germany}, \texttt{email:linde@mathematik.uni-muenchen.de}; \newline current address: \textsl{Department of Pure Mathematics, University of Liverpool, Peach Street, Liverpool. L69 7ZL, United Kingdom}, \texttt{email:scsa0296@liverpool.ac.uk};\newline supported by the Marie Curie Fellowship ref. HPMT-CT-2001-00277.}  \date{} }
\begin{document}
\maketitle
\setcounter{secnumdepth}{4}
\setcounter{section}{0}
\begin{abstract}
It is shown that a certain representation of the Heisenberg type Krichever-Novikov algebra gives rise to a state field correspondence that is quite similar to the vertex algebra structure of the usual Heisenberg algebra. Finally a definition of Krichever-Novikov type vertex algebras is proposed and its relation to the ''classical'' vertex algebra is discussed.
\end{abstract}
\section{Introduction}
Vertex algebras can be considered as an algebraic counterpart of two-dimensional conformal field theory \cite{kac}. Typical examples of vertex algebras are the Heisenberg vertex algebra, lattice vertex algebras and vertex algebras of affine Kac-Moody algebras. They correspond to toroidal compactification, and  WZNW-models, respectively.
\\ There are some different approaches to extend conformal field theory on Riemann surfaces, as for example \cite{knty} (see also \cite{por1}\cite{por2} for further developments) and \cite{tuy}. From the point of view of vertex operators see \cite{fre}.
\\ A vertex algebra $V$ ($V$ a vector space) consists in principle of a state field correspondence $V\rightarrow End(V)[[z,z^{-1}]]$ that maps elements of $a\in V$ to fields $\sum_{n}a_{n}z^{-n-1}$ where $a_{n}\in End(V)$, and $a_{n}v=0$ for $n>>0$.
\\ Krichever and Novikov generalized in \cite{kn1}\cite{kn2}\cite{kn3} the notion of affine algebras in the following way:
\\ If we identify the associative algebra $\C[t,t^{-1}]$ of Laurent polynomials with the algebra consisting of meromorphic functions on the Riemann sphere $\hat{\C}$ which are holomorphic outside $0$ and $\infty$ ($t$ is a quasi-global coordinate on $\hat{\C}$), then we can generalize the current algebra to higher genus Riemann surfaces where $\C[t,t^{-1}]$ is replaced by the associative algebra $\mathcal{A}(X,P_{\pm})$ of meromorphic functions that are holomorphic outside two distinguished points. From this point of view the central extended current algebra
\[ \hat{\la{g}} = \la{g}\otimes \C[t,t^{-1}] \oplus \C K\]
is replaced by
\[ \hat{\la{g}} = \la{g}\otimes \mathcal{A} \oplus \C K. \]
Krichever-Novikov algebras have numerous relations  to the fundamental problems of geometry, analysis and mathematical physics.
\\ We only mention very briefly their relationship to moduli of stable vector bundles, the geometric Langlands correspondence \cite{sh2} and to Knizhnik-Zamolodchikov equations \cite{schsh1}.
\\ Applying the considerations of constructing Krichever-Novikov type algebras to the definition of vertex algebras raises the question of whether it is possible to generalize vertex algebras to higher genus Riemann surfaces. Or in other words, if the notion of vertex algebras is a generalization of the notion of a unital commutative associative algebra where the multiplication depends on a parameter \cite{bakkac} then in a higher genus vertex algebra the complex parameter should encode the geometry of the corresponding Riemann surface.
\\ The fields in vertex algebras of the form
\[ \sum_{n} a_{n}z^{-n-\Delta_{a}} \]
must then be replaced by fields of the form
\[ \sum_{n} a_{n} f_{\Delta_{a}}^{n}(P) \]
where $f_{\Delta_{a}}^{n}(P)$ is a section of the $\Delta_{a}$ times tensorized canonical bundle $K$, and with possible poles in the points $P_{+}$ and $P_{-}$.
\\ That means that the state field correspondence should now be given by a linear map $V\rightarrow End(V)[[f_{\lambda}^{n}]]_{\lambda\in\Z_{\geq 0},n\in\Z}$ ($V$ a vector space).
\\ The state field correspondence for vertex algebras has to satisfy certain properties such as the vacuum axiom, the translation axiom and the locality axiom. These properties have to be reformulated for a higher genus vertex algebra. The operator product has to be reformulated in a suitable way as well.
For example for generating fields of two elements of an affine algebra we get
\[ [J^{a}(z),J^{b}(w)] = [J^{a},J^{b}](w)\delta(z-w) + (J^{a},J^{b})\partial_{w} \delta(z-w) \]
This should be replaced by
\[ [a(P),b(Q)] = [a,b](P)\Delta(P,Q) + (a,b) d_{P}\Delta(P,Q) \]
where $\Delta(P,Q)$ means the delta distribution for higher genus. It is defined as the formal sum of the expansions of the Szeg\"{o} kernel.
\\ This paper is organized as follows:
\\ In section 2 we introduce some notations for Krichever-Novikov algebras. We use the old-fashioned version of bases. This means the bases are in some sense more symmetric than e.g. in \cite{sh3}. We also discuss the physical meaning of level lines as it was proposed by Krichever and Novikov in \cite{kn2}. Namely the level lines can be considered as the snapshots of propagating strings. 
\\ One of the key features in Conformal field theory on compact Riemann surfaces is the Szeg\"{o} Kernel. The correlation functions of fields on Riemann surfaces are expressions in terms of the Szeg\"{o} Kernel \cite{raina1}\cite{knty}\cite{fay2}\cite{vv}. We present for our purposes an expansion of a certain Szeg\"{o} kernel.
\\ In section 3 the definitions of admissible representations and of fields are given, and derivatives of fields  are discussed. For a geometric approach to the representation theory of affine Krichever-Novikov algebras by so called framed bundles see \cite{sh2}. 
\\ In section 4 we define the normal ordered product of fields on a Riemann surface by mimicking the normal ordered product for the formal calculus (see \cite{kac}). From this point of view it is very natural to define the normal ordered product on a compact Riemann surface as a formal expression where the Szeg\"{o} kernel is involved.
\\ In section 5 we define a state-field correspondence on the space $V$ of a suitable representation of Heisenberg type KN-algebra, i.e. a linear map that maps each element of $V$ to a formal series of endomorphisms of $V$.
\\ In section 6 some properties of this map are discussed, namely the vacuum property and the translation property.
\\ In section 7 the analogous locality axiom property is discussed. In the vertex algebra case locality means that two fields satisfy 
\[ [a(z),b(w)](z-w)^{N} = 0 \]
As a suitable analogue in the higher genus case the locality condition is supposed to be replaced by the condition that the product with the Schottky-Klein prime form
\[ [a(P),b(Q)](E(P,Q))^{N} \]
gives a multivalued but regular section on $X\times X$.
\\ In section 8 a definition of vertex algebras for Krichever-Novikov algebras is suggested that might give the appropriate framework for considering conformal field theory related questions. It is also shown that the given definition coincides with the definition of a ''classical'' vertex algebra in the case $g=0$. The notion of a vertex algebra is generalized to chiral algebras in the sense of Beilinson and Drinfeld \cite{bd}. Some possible connections of our definition with this notion are discussed.
\\ In section 9 obvious extensions to affine KN-algebras are discussed.
\\ In a forthcoming paper the material will be carried out in more detail.
\section{Notations}
\subsection{Krichever-Novikov Bases}
Let $X$ be a compact Riemann surface of genus $g\geq 1$ and $P_{\pm}$ two distinguished points in general position. Denote by $K$ the canonical line bundle on $X$.
\\ Let be $\lambda\in\Z$. Denote by $\mathcal{F}^{\lambda}$ the infinite dimensional vector space of global meromorphic sections of $K^{\lambda}:= K^{\otimes \lambda}$ which are holomorphic outside $P_{\pm}$. 
\\ For $\lambda=0$, i.e. for the space of meromorphic functions holomorphic outside $P_{\pm}$ we will write $\mathcal{A}:=\mathcal{A}(X,P_{\pm}):=\mathcal{F}^{0}$. By multiplying sections in $\mathcal{F}^{\lambda}$ by functions we again obtain sections in $\mathcal{F}^{\lambda}$. In this way the space $\mathcal{A}$ becomes an associative algebra and the $\mathcal{F}^{\lambda}$ become modules over $\mathcal{A}$.
\\ The elements in $\mathcal{F}^{-1}$ that are vector fields operate on the spaces $\mathcal{F}^{\lambda}$ by taking the Lie derivative. In local coordinates
\begin{equation}
\label{eq:der}
\nabla_{e}(g)|  :=  \left( e(z)\frac{d}{dz} \right) .\left( g(z)(dz)^{\lambda}\right) = \left( e(z)\frac{dg}{dz}(z)+ \lambda g(z)\frac{de}{dz}(z) \right) (dz)^{\lambda}.
\end{equation}
Here $e\in\mathcal{F}^{-1}$ and $g\in\mathcal{F}^{\lambda}$. The space $\mathcal{F}^{-1}$ becomes a Lie algebra with respect to the above equation, and the spaces $\mathcal{F}^{\lambda}$ become Lie modules \cite{schsh1}.
\begin{dfn}[Krichever-Novikov pairing]
The \textbf{Krichever-Novikov pairing} is the pairing between $\mathcal{F}^{\lambda}$ and $\mathcal{F}^{1-\lambda}$ given by
\[ \mathcal{F}^{\lambda}\times \mathcal{F}^{1-\lambda}\rightarrow \C, \hspace{0.4cm} (f,g):= \frac{1}{2\pi i}\int_{C} f\cdot g, \]
where $C$ is an arbitrary non-singular level line. (see section \ref{sec:level} for the definition)
\end{dfn}
We consider now certain bases of the vector spaces $\mathcal{F}^{\lambda}$, called \textbf{Krichever-Novikov bases} or KN-bases for short. 
\\ Let be $g>1$. The KN-bases are supposed to have the following behaviour around the points $P_{\pm}$:
\begin{equation}
\label{eq:knbase}
f_{\lambda,n}(z_{\pm})=\alpha_{n,\pm}^{\lambda} z_{\pm}^{\pm n-s_{\lambda}}\left( 1+O(z_{\pm})\right) (dz_{\pm})^{\lambda} 
\end{equation}
where 
\begin{equation} s_{\lambda}= \frac{(1-2\lambda)g}{2}+\lambda \end{equation}
with constants $\alpha_{\pm}^{\lambda}\in\C$; $n\in\Z$, if $g$ even, and $n\in\Z+\frac{1}{2}$, if $g$ odd.
\\ By these conditions the $f_{\lambda, n}$ are uniquely determined.
\\ For the case $\lambda=0$ we have some modifications: 
\\ Denote by $A_{n}(P)\in\mathcal{A}$ the functions with the local behaviour:
\begin{eqnarray}
A_{n}(z_{\pm}) & = & \alpha_{n,\pm}^{0} z_{\pm}^{\pm n-\frac{g}{2}}\left( 1+O(z_{\pm})\right)  \mbox{ if } |n|\geq \frac{g}{2}+1\\
A_{n}(z_{\pm}) & = & \alpha_{n,\pm}^{0} z_{\pm}^{\pm n-\frac{g}{2}\left\{ +0/-1\right.} \left( 1+O(z_{\pm})\right) \mbox{ if } n\in\{ -\frac{g}{2},...,\frac{g}{2}-1\} \\
\label{eq:a1}
A_{\frac{g}{2}} & = & 1
\end{eqnarray}
We set $\alpha_{n,+}^{0}=1$ for all $n$.
\\ For the case $\lambda=1$ we have some modifications as well: 
\\ Denote by $\omega^{n}(P)\in\mathcal{F}^{1}$ the 1-forms with the local behaviour:
\begin{eqnarray}
\omega^{n}(z_{\pm}) & = & \alpha_{\pm}^{1} z_{\pm}^{\mp n+\frac{g}{2}-1}\left( 1+O(z_{\pm})\right)dz  \mbox{ if } |n|\geq \frac{g}{2}+1\\
\omega^{n}(z_{\pm}) & = & \alpha_{\pm}^{1} z_{\pm}^{\mp n+\frac{g}{2}\left\{ -1/+0\right.} \left( 1+O(z_{\pm})\right)dz \mbox{ if } n\in\{ -\frac{g}{2},...,\frac{g}{2}-1\} \\
\omega^{\frac{g}{2}} & : & \mbox{defining differential of the level lines.}
\end{eqnarray}
For $g=1$ we have a holomorphic 1-form $\omega$ (without zeros and poles) and we can write $f_{\lambda,n}= A_{n}\cdot (\omega)^{\lambda}$ (where $A_{n}$ is defined as above) (see e.g. \cite{bo} for details).
\begin{dfn}[Residue]
Let be $f(P)$ a 1-form regular outside the points $P_{\pm}$. Then we define
\begin{equation}
\mbox{Res }\left( f(P)\right) = \frac{1}{2\pi i} \int_{C_{s}} f(P) 
\end{equation}
where $C_{s}$ is a cycle cohomologous to a small circle around the point $P_{+}$.
\end{dfn}
Let be the constants $\alpha_{\pm}^{\lambda}$ defined such that we have the duality relations:
\[ \mbox{Res }\left( f_{\lambda,n}f_{-\lambda+1}^{m} \right) = \delta_{n,m} \]
where $f_{\lambda}^{n}\stackrel{def}{=}f_{\lambda,-n}$.
\begin{lem} The algebra $\mathcal{A}$ acts on the bases of $\mathcal{F}^{\lambda}$ as follows:
\begin{equation} \label{eq:beta} A_{n}(P)f_{\lambda}^{m}(P) = \sum_{k}\beta_{nk}^{\lambda, m} f_{\lambda}^{k}(P) \end{equation}
with
\[ \beta_{nk}^{\lambda, m} \stackrel{def}{=} \mbox{Res }\left( A_{n}(P)f_{\lambda}^{m}(P)f_{1-\lambda,k}(P) \right) \]
and $\beta_{nk}^{\lambda, m}\neq 0$ only for finitely many $k$. More precisely there exist constants $c_{1},c_{2}$ such that $\beta_{nk}^{\lambda, m}\neq 0 \Rightarrow n-m-c_{1}\leq k\leq n-m+c_{2}$.
\end{lem}
We obtain the observation:
\begin{cor}
\label{res:alpha}  The $\beta$s are related in the following way:
\begin{equation} \beta_{nk}^{\lambda,m} = \beta_{n,-m}^{1-\lambda, -k}. \end{equation}
For $\lambda=1$ with $\alpha_{nk}^{m} := \beta_{nk}^{1, m}$ we have:
\begin{eqnarray}
A_{n}(P)\omega^{m}(P) & = & \sum_{k}\alpha_{nk}^{m} \omega^{k}(P) \\
A_{n}(P)A_{k}(P) & = & \sum_{m}\alpha_{nk}^{m}A_{m}(P)
\end{eqnarray}
\end{cor}
\subsection{Level Lines}
\label{sec:level}
Let $X$ be a compact Riemann surface of genus $g\geq 1$ as above. Let be $\alpha_{1},...\alpha_{g},\beta_{1},..,\beta_{g}$ a canonical homotopy basis for $X$. $\alpha_{1},...\alpha_{g},\beta_{1},..,\beta_{g}$ determines a basis $\{ v_{1},..,v_{g}\}$ for the holomorphic differentials on $X$ with normalized periods
\[ \left(\int_{\alpha_{j}}v_{k}, \int_{\beta_{j}}v_{k}\right) = \left(E,\Omega \right) \]
for some $\Omega\in\C^{g\times g}$, with $Im\Omega$ positive definite (E is the identity $g\times g$-matrix).
\\ We have the following well known result (e.g. \cite{fay2}):
\begin{prop}
\label{sec:level}
There is a unique differential $\rho$ with $\mbox{Res}_{P_{\pm}}\rho=\pm 1$ and with poles of order $1$ in $P_{\pm}$, and with purely imaginary periods. This differential can be given in analytic terms by
\begin{equation} \omega_{P_{\pm}}(P) = d\left(log\frac{E(P,P_{+})}{E(P,P_{-})}\right) - 2\pi i\sum_{ij=1}^{g} Im\left( \int_{P_{-}}^{P_{+}} v_{i}\right) \left( Im\Omega \right)_{ij}^{-1} v_{j}(P) \end{equation}
where $E(P,Q)$ is the Schottky-Klein prime form.
\end{prop} 
According to the differential $\rho$ we can define the level lines \cite{kn1}:
\begin{dfn}[Level Lines] Let be $P_{0}$ a point on $X$ different from $P_{\pm}$ Then the level line of level $\tau\in\R$ is defined by
\begin{equation} C_{\tau} =\left\{ P\in X : \mathfrak{Re}\int_{P_{0}}^{P} \rho \right\} \end{equation}
\end{dfn}
The level lines can be considered as snapshots of propagating strings (\cite{kn2}). With this interpretation we get a sort of time axis on the Riemann surface. For $P\rightarrow P_{-}$ the level lines become circles around $P_{-}$ and for $P\rightarrow P_{+}$ they become circles around $P_{+}$.
\\ We can distinguish points on the Riemann surface by the following criterion: $\tau(P)\geq\tau(Q)$ $(\tau(P)>\tau(Q))$ means that the level line of $P$ has a lower or equal (bigger) level than the level line corresponding to $Q$.
\\ For $g=0$ we can identify $P_{+}$ with zero and $P_{-}$ with $\infty$, then the above defined $\rho$ corresponds to $\frac{1}{z}dz$ on the complex plane and we get as level lines on $\C$ circles around zero given by $\mathfrak{Re}\int_{1}^{z}\frac{1}{z}dz$.
%
\subsection{Szeg\"{o} Kernel}
\begin{dfn}[Szeg\"{o} Kernel]
Let $L\in Pic^{g-1}$ be a line bundle on $X$ of degree $g-1$ and $H^{0}(X,L)\neq 0$.
The Szeg\"{o} Kernel is defined as the unique section (up to a constant):
\[ \mathfrak{s}_{L}\in H^{0}(X\times X, L\sqtimes K\otimes L^{*}(\Delta)) \]
\end{dfn}
The fact that the Szeg\"{o} kernel is unique is well known (see e.g. \cite{raina1}).
\\ In analytic terms we can write
\[ \mathfrak{s}_{L}(P,Q) = \frac{\theta[L](\int_{Q}^{P}v)}{\theta[L](0)E(P,Q)} \]
Consider now the following (convergent) expansions that can be considered as expansions of a certain Szeg\"{o} kernel:
\begin{eqnarray}\label{eq:szegoexp}
i_{P,Q}S(P,Q) & = & \sum_{n=\frac{g}{2}}^{\infty} A_{n}(P) \omega^{n}(Q) \\
i_{Q,P}S(P,Q) & = & -\sum_{n=-\infty}^{\frac{g}{2}-1} A_{n}(P) \omega^{n}(Q) 
\end{eqnarray}
where $i_{P,Q}$ $(i_{Q,P})$ means the expansion in the region $\tau(P)<\tau(Q)$ $(\tau(P)>\tau(Q))$ (see \cite{kn2}).
%
\section{Representations and Fields}
\subsection{The Heisenberg Algebra of Krichever-Novikov Type}
Let $\la{h}$ be an $l$-dimensional complex vector space, let $\la{h}$ be equipped with a non-degenerate, symmetric, bilinear form $(\cdot |\cdot)$ and an orthonormal basis $a^{1},...,a^{l}$.
\\ The Heisenberg algebra is the affinization of $\la{h}$:
\[ \hat{\la{h}} = \la{h}\otimes\C[t,t^{-1}]\oplus \C K\]
with brackets $[K,\hat{\la{h}}] =0$ and
\[ [a_{n},b_{m}] = (a|b)n\delta_{n,-m} K \hspace{0.5cm} \mbox{where }a_{n}=a\otimes t^{n}. \]
%
We can proceed in the KN-case (see e.g. \cite{schsh1}) especially for the case $\la{h}=\C$:
\\ The \textbf{Heisenberg algebra of KN-type} is defined by:
\[ \hat{\mathcal{A}} := \hat{\mathcal{A}}(X,P_{\pm}) : = \C\otimes \mathcal{A}\oplus \C K \]
Define by $a_{n}:= a\otimes A_{n}$ the elements of $\C\otimes \mathcal{A}$. Then we have the brackets:
\[ [a_{n},a_{m}] = \gamma_{nm} K \hspace{0.5cm} [a_{n},K] = 0 \]
where
\[ \gamma_{nm} = \mbox{Res }\left( A_{n} dA_{m} \right) \]
and
\begin{equation}\label{eq:gamma}  \begin{array}{ccc} \gamma_{nm}=0& \mbox{for}& |n+m|>g, |n|,|m|>\frac{g}{2} \\ \gamma_{nm}=0& \mbox{for}& |n+m|>g+1, |n|\mbox{ or }|m|\leq\frac{g}{2} \end{array} \end{equation}
Let be $\hat{\mathcal{A}}_{n}=\C\otimes A_{n}$ for all $n\neq \frac{g}{2}$, and  $\hat{\mathcal{A}}_{\frac{g}{2}}=\C\oplus\C K$. Then we can decompose
\[ \hat{\mathcal{A}} = \bigoplus_{n=\frac{g}{2}\mod \Z} \hat{\mathcal{A}}_{n}. \]
Denote by $\mathcal{A}_{+}$ $(\mathcal{A}_{-})$ the functions that have zeros at the points $P_{+}$ $(P_{-})$. The functions $A_{n}$ with $n>\frac{g}{2}$ $(n<-\frac{g}{2})$ form a basis in $\mathcal{A}_{+}$ $(\mathcal{A}_{-})$. Denote by $A_{0}$ the space of functions $A_{n}$ with $|n|\leq\frac{g}{2}$. We get therefore the decomposition:
\[ \mathcal{A} = \mathcal{A}_{+} \oplus \mathcal{A}_{0} \oplus \mathcal{A}_{-}. \]
\subsection{Representations}
\begin{dfn}[Admissible Representation] Let $V$ be a vector space. An \textbf{Admissible Representation} of the Heisenberg Algebra $\hat{\mathcal{A}}(X,P_{\pm})$ is a linear map 
\[ \pi: \hat{\mathcal{A}}(X,P_{\pm}) \rightarrow End (V) \]
where for all $v\in V$ exists an $n_{0}\in\Z$ such that $a_{n}v=0$ for all $n\geq n_{0}$.
\end{dfn}
Example: Fock representation of the Heisenberg algebra.
\\ For an element $\eta$ in the dual space $\la{h}^{*}$ we consider the Fock space defined by the induced module
\[ V_{\eta} = U(\hat{\la{h}})\otimes_{U(\la{h}\otimes \C[t]\oplus\C K)} \C, \]
where $\C$ is the one-dimensional space annihilated by $\la{h}\otimes t\C[t]$ and on which $K$ acts as the identity and $\la{h}\otimes t^{0}$ via the character $\eta$. $U$ denotes the universal enveloping algebra. 
\\ Especially for $l=1$ the space $V_{0}$ is spanned by the vectors
\[ a_{-n_{1}-1}...a_{-n_{M}-1}v_{0} \mbox{ where } n_{i}>0. \]
\subsection{Fields}
The usual definition of fields is as follows (see \cite{kac}):
\begin{dfn} Let $V$ be a vector space. A formal series $a(z)\in End(V)[[z,z,^{-1}]$ is said to be a \textbf{field}, if
\[ \mbox{for all }v\in V\exists n_{0}: a_{n}v_{0} =0 \mbox{ for all }n\geq n_{0} \]
in other words: the series $\sum_{n}(a_{n}v)z^{-n-1}\in V[[z]][z^{-1}]$ is a formal Laurent series.
\end{dfn}
We extend this definition to the higher genus case.
\begin{dfn}[Field] A \textbf{field} of weight $\lambda$ is a formal series $a(P)\in End(V)[[(f_{\lambda}^{n}(P))_{n\in\Z}]]$,
\[ a(P) = \sum_{n} a_{n} f_{\lambda}^{n}(P) \]
with the property:
\[ \forall v\in V\exists n_{0}\forall n\geq n_{0}: a_{n}v = 0. \]
\end{dfn}
\subsection{Derivatives of Fields}
For fields $a(z)\in End(V)[[z,z,^{-1}]]$ with
\[ a(z) = \sum_{n}a_{n} z^{-n-1} \]
we have as  formal derivative
\[ \partial a(z) = \sum_{n} (-n-1)a_{n} z^{-n-2} = \sum_{n} (-n) a_{n-1} z^{-n-1}, \]
and iterated derivation gives
\[ \frac{1}{k!} \partial^{k} a(z) = \sum_{n} (-1)^{k}{n\choose k} a_{n-k} z^{-n-1} = \sum_{n} a^{(k)}_{n} z^{-n-1}, \]
where $a^{(k)}_{n} = (-1)^{k}{n\choose k}a_{n-k}$. And because of the properties of the binomial coefficients we get
\begin{equation}\label{eq:zerox} a^{(k)}_{n}=0 \mbox{ for } k>n\geq 0. \end{equation}
For the KN-fields the derivation of a field of conformal dimension $\lambda=1$ is defined by
\begin{equation} \nabla_{e}a(P) = \sum_{n}a_{n}(\nabla_{e} \omega^{n})(P) \end{equation}
where $e$ is a vector field regular outside $P_{\pm}$ (see eq. (\ref{eq:der})).
\\ The basis elements $e_{k}$ look according to eq. (\ref{eq:knbase}) locally around $P=P_{+}$:
\[ e_{n}(z_{+}) = \alpha^{-1}_{n,+} z_{+}^{n-\frac{3g}{2}+1}\left(1+O(z_{+}) \right) \frac{d}{dz} \]
\begin{dfn} We consider especially the case $n=\frac{3g}{2}-1$: 
\[ \nabla \stackrel{def}{=} \nabla_{e_{\frac{3g}{2}-1}}. \]
\end{dfn}
\vspace{0.2cm}
We get according to this definition:
\[ \nabla a(P) = \sum_{n,u} \zeta_{u}^{n}a_{n}\omega^{u}(P) \]
where the coefficients $\zeta_{u}^{n}$ are given by
\[ \zeta_{u}^{n} = \mbox{Res }\left( \omega^{n}(Q) e_{\frac{3g}{2}-1}(Q) dA_{u}(Q) \right). \]
\begin{prop} Because $A_{\frac{g}{2}} = 1$ (eq. (\ref{eq:a1})) we get especially:
\begin{equation} \label{eq:zeta0} \zeta_{\frac{g}{2}}^{n} = 0 \forall n.  \end{equation}
There exists $C>0$ such that
\begin{equation}\label{eq:zetarel} \zeta_{n}^{m}\neq 0 \Rightarrow n-1\leq m \leq C+n \end{equation}
\end{prop}
Iterated derivation gives:
\[ \nabla^{k} a(P) = \sum_{u}\sum_{u_{k-1},...,u_{1},n} \zeta_{u}^{u_{k-1}}...\zeta_{u_{1}}^{n} a_{n} \omega^{u}(P)  = \sum_{u} a^{(k)}_{u}\omega^{u}(P) \]
i.e.
\begin{equation} a^{(k)}_{u} \stackrel{def}{=} \sum_{u_{k-1},...,u_{1},n} \zeta_{u}^{u_{k-1}}...\zeta_{u_{1}}^{n} a_{n} \end{equation}
These sums are finite due to eq. (\ref{eq:zetarel}).
\\ We write $a^{(k)}_{u}$ in the form
\[ a^{(k)}_{u} = \sum_{j} q_{u}^{(k),j}a_{j} \]
Because of eq (\ref{eq:zeta0}) we get a similar result as eq. (\ref{eq:zerox}):
\begin{equation}\label{eq:zero} q_{u}^{(k),j} = 0 \mbox{ for } j<\frac{g}{2}, \frac{g}{2}\leq u\leq \frac{g}{2}+k-1 \end{equation}
%
%
%
\section{Normal Ordered Product}
For the case of fields in the usual sense the normal ordered product is defined by (see \cite{kac}):
\[ :a(w)b(w): = \mbox{Res} \left(a(z)b(w)i_{z,w} \frac{1}{z-w} +  b(w)a(z)i_{w,z} \frac{1}{z-w} \right) \]
We get for the coefficients of two fields $a(z),b(z)$ where $a(z)$ is of conformal dimension $1$, and $b(z)$ is of conformal dimension $\lambda$ for the coefficients $:a(w)b(w): = \sum_{n} :a(w)b(w):_{n} w^{-n-\lambda-1}$:
\begin{equation}\label{eq:g0nop} :a(w)b(w):_{n} = \sum_{j<0}a_{j}b_{n-j} + \sum_{j\geq 0} b_{n-j}a_{j}. \end{equation}
This motivates the following definition for higher genus:
\begin{dfn}[Normal Ordered Product] \label{dfn:nop} Let be $a(P)\in End V[[(\omega^{n})_{n\in\Z}]]$ a field of weight $1$ and let be $b(P)\in End V[[(f_{\lambda}^{n})_{n\in\Z}]]$ a field of weight $\lambda$.
\\ The \textbf{Normal Ordered Product} is defined by
\begin{equation} :a(Q)b(Q): = \mbox{Res} \left(a(P)b(Q)i_{Q,P} S(P,Q) +  a(Q)b(P)i_{P,Q} S(P,Q) \right) \end{equation}
where the symbol $i_{P,Q}S(P,Q)$ means expansion as in eq. (\ref{eq:szegoexp}).
\end{dfn}
\begin{prop} 
\label{prp:nop} The coefficients of $:a(Q)b(Q): = \sum_{n} :a(Q)b(Q):_{n}f_{\lambda+1}^{n}(P)$ can be written as
\begin{equation} :a(Q)b(Q):_{n} = \sum_{m, j<\frac{g}{2}}a_{j}b_{m}l^{jm}_{n} + \sum_{m, j\geq \frac{g}{2}}b_{m}a_{j}l^{jm}_{n} \end{equation}
where
\[ l^{jm}_{n} \stackrel{def}{=} l^{jm}_{n,(\lambda)}  \stackrel{def}{=} \mbox{Res }\left(\omega^{j}(P) f_{\lambda}^{m}(P) f_{-\lambda}^{-n}(P)\right)  \]
and we have the relation
\begin{equation}\label{eq:ljmn} l^{jm}_{n} \neq 0 \Rightarrow n-j-\frac{g}{2}\leq m\leq n-j+\frac{g}{2}. \end{equation}
\end{prop}
In the case $g=0$ we have the relation
\[ l^{jm}_{n} = \mbox{Res }\left( z^{-j-1}z^{-m-\lambda}z^{n+\lambda}dz \right) = \delta_{n-j,m}, \]
so that we obtain eq. (\ref{eq:g0nop}).
\\ \textsl{Proof of the proposition.}
We use the expansion of the Szeg\"{o}-Kernel and duality and get:
\[ :a(Q)b(Q):  = \]
\[ = \sum_{m, j<\frac{g}{2}}a_{j}b_{m}f_{\lambda}^{m}(Q)\omega^{j}(Q) + \sum_{m, j\geq \frac{g}{2}}b_{m}a_{j}f_{\lambda}^{m}(Q)\omega^{j}(Q) \]
\[ \stackrel{(*)}{=} \sum_{m, j<\frac{g}{2}}a_{j}b_{m}\sum_{n}l^{jm}_{n} f_{\lambda+1}^{n}(Q) + \sum_{m, j\geq \frac{g}{2}}b_{m}a_{j}\sum_{n}l^{jm}_{n} f_{\lambda+1}^{n}(Q) \]
\[ = \sum_{n}\left( \sum_{m, j<\frac{g}{2}}a_{j}b_{m}l^{jm}_{n} + \sum_{m, j\geq \frac{g}{2}}b_{m}a_{j}l^{jm}_{n}  \right)  f_{\lambda+1}^{n}(Q) \]
ad (*): We used the relation:
\[ \omega^{j}(Q)f_{\lambda}^{m}(Q) = \sum_{n} l^{jm}_{n} f_{\lambda+1}^{n}(Q). \hspace{1cm}\square\] 
\\ There are some consequences of this definition of the normal ordered product.
\begin{thm} The normal ordered product of a field of weight $1$ and a field of weight $\lambda$ is a field of weight $\lambda+1$.
\end{thm}
Furthermore we can define an iterated normal ordered product:
\\ Let be $a^{1}(P),...,a^{N}(P)\in End V[[(\omega^{n})_{n\in\Z}]]$ fields of weight $1$, and let be $b(P)\in End V[[(f_{\lambda}^{n})_{n\in\Z}]]$ a field of weight $\lambda$. 
\\ Then the iterated normal ordered product (from the left) is
\begin{equation} :a^{1}(P):a^{2}(P):...:a^{N}(P):b(P):...: \end{equation}
and it is (due to the above theorem) a field of weight $\lambda+N$.
%
%
\section{State Field Correspondence}
\subsection{Vertex Algebras}
A \textbf{vertex algebra} is a Triple $(V,Y,v_{0})$ consisting of a $\Z_{\geq 0}$-graded vector space $V=\bigoplus_{n\in\Z_{\geq 0}} V_{n}$, a linear map $Y:V\rightarrow End(V)[[z,z,^{-1}]]$, and a distinguished vector $v_{0}\in V_{0}$. The statefield correspondence for a vector $A\in V$ of \textbf{conformal dimension} $\Delta_{A}$ is as follows:
\[ A\mapsto Y(A,z) = \sum_{n}A_{n} z^{-n-\Delta_{A}}. \]
We have the vacuum axiom: $Y(v,z)v_{0}|_{z=0} = v_{-\Delta_{A}}$, in other words $A_{n}v_{0}=0$ for all $n > -\Delta_{A}$, and $A_{-\Delta_{A}}v_{0} = 0$.
\\ We have the translation axiom: $Tv_{0} = 0$, and $\partial Y(A,z) = [T,Y(A,z)]$.
\\ And finally the locality axiom: There exists an $N$ such that
\[ [Y(A,z), Y(B,w)](z-w)^{N} = 0. \]
\textsl{Example: The Heisenberg vertex algebra}
\\ The state field correspondence for the above given representation of the Heisenberg algebra
\begin{eqnarray*}
a_{-n_{1}}...a_{-n_{M}} & \mapsto & Y( a_{-n_{1}}...a_{n_{M}} ,z) = \\
& = & :\partial^{(n_{1}-1)}a(z)... \partial^{(n_{M}-1)}a(z): 
\end{eqnarray*}
defines a vertex algebra structure (see e.g. \cite{fre}).
\subsection{Higher Genus State Field Correspondence}
Let be $\pi: \hat{\mathcal{A}}(X,P_{\pm})\rightarrow End(V)$ an admissible representation with the following properties:
\begin{enumerate}
\item There is a vector $v_{0}\in V$ such that $a_{\frac{g}{2}+h}v_{0} = 0 $ $\forall h\in\Z_{\geq 0}$. $K$ acts as the identity on $V$.
\item $V$ is spanned by the vectors 
\begin{equation}\label{eq:span} a_{-n_{1}+\frac{g}{2}}a_{-n_{2}+\frac{g}{2}}....a_{-n_{M}+\frac{g}{2}} v_{0} \end{equation}
where $n_{1}\geq ...\geq n_{M}> 0$.
\item Gradation: Let $v\in V$ be a vector as in eq. (\ref{eq:span}). The degree of $v$ is defined by
\[ \deg(v) = \sum_{j=1}^{M} n_{j}. \]
\end{enumerate}
\begin{dfn}
Let be given a representation $\pi: \hat{\mathcal{A}}(X,P_{\pm})\rightarrow End(V)$ as above.
\\ Define a linear map (the state-field correspondence) 
\[ \mathcal{Y}: V\rightarrow End(V)[[(f_{\lambda}^{n})_{n\in\Z, \lambda\in\Z_{>0}} ]] \]
by
\begin{eqnarray*}
a_{-n_{1}-\frac{g}{2}-1} ... a_{-n_{M}-\frac{g}{2}-1}v_{0} & \mapsto & \mathcal{Y}(a_{-n_{1}-\frac{g}{2}-1} ... a_{-n_{M}-\frac{g}{2}-1}v_{0}, P) \\
& \stackrel{def}{=} & :\nabla^{n_{1}}a(P)... \nabla^{n_{M}}a(P): 
\end{eqnarray*} 
\end{dfn}
In the next two sections we are going to show that this linear map satisfies some properties that are analogous to the usual vertex algebra axioms.
%
%
\section{Vacuum Property}
Similar to the Heisenberg vertex algebra we have the following result:
\begin{thm}
The Fields
\[ :\nabla^{n_{1}-1}a(P) \nabla^{n_{2}-1}a(P) .... \nabla^{n_{M}-1}a(P): \]
satisfy the property:
\begin{equation} :\nabla^{n_{1}-1}a(P) \nabla^{n_{2}-1}a(P) .... \nabla^{n_{M}-1}a(P):_{n}v_{0} = 0 \mbox{ for } n> -s_{M}  \end{equation}
and
\begin{eqnarray*}
 :\nabla^{n_{1}-1}a(P) .... \nabla^{n_{M}-1}a(P):v_{0}|_{P=P_{+}} & = & :\nabla^{n_{1}-1}a(P) .... \nabla^{n_{M}-1}a(P):_{-s_{M}}v_{0} \\
& = & C\cdot \left( a_{-n_{1}+\frac{g}{2}}....a_{-n_{M}+\frac{g}{2}}\right) v_{0} + ...
\end{eqnarray*}
where $...$ means lower degree vectors. $C$ is a scalar.
\end{thm}
\begin{lem} For the field $a(P):= \mathcal{Y}(a_{\frac{g}{2}-1}v_{0},P)=\sum_{n}a_{n}\omega^{n}(P)$ we have
\label{res:a0}
\[ a(P)v_{0}|_{P=P_{+}} = c\cdot a_{\frac{g}{2}-1}v_{0} \]
where $c\in\C$.
\end{lem}
Proof. 
\[ a(P)v_{0} = \sum_{n} a_{n}v_{0} \omega^{n}(P) = \sum_{n<\frac{g}{2}} a_{n}v_{0} \omega^{n}(P) \]
because $a_{\frac{g}{2}+x}v_{0} = 0$ for all $x=0,1,2,...$.
\\ The $\omega^{n}(P)$ look locally around $P_{+}$: $z_{+}^{-n+\frac{g}{2}-1}(1+O(z_{+}))dz_{+}$.
\\ For $n\leq \frac{g}{2}-1$ we get non-negative powers and therefore
\begin{equation}
\label{eq:omega}
\omega^{n}(P_{+}) =\left\{ \begin{array}{cc} \alpha^{\frac{g}{2}-1}_{1,+} & \mbox{for } n=\frac{g}{2}-1 \\ 0 & \mbox{for } n<\frac{g}{2}-1 \end{array}\right. 
 \end{equation} 
%
\\ For the derivations of the field $a(P)$ we get the following result:
\begin{prop}
\label{prp:der}
\[ \nabla^{k}a(P)v_{0}|_{P=P_{+}} = \eta\cdot a_{-k+\frac{g}{2}-1}v_{0} + .... \]
where $...$ means summands of lower degree $(\eta\in\C)$.
\end{prop}
This is a consequence of eq. (\ref{eq:zero}).  \hfill{$\square$}
\\ We turn now to the proof of the theorem:
Let be $M=1$ and $n_{M}=n_{1}=1$. Then we have due to (Lemma \ref{res:a0}):
\[ a(P) v_{0}|_{P=P_{+}} = a_{\frac{g}{2}-1}v_{0} \] 
For $n_{M}=n_{1}>1$ we have due to Propos. \ref{prp:der}:
\[ \nabla^{k}a(P)v_{0}|_{P=P_{+}} = \eta\cdot a_{-k+\frac{g}{2}-1}v_{0} + .... \]
and we obtain the assertion of the theorem because $s_{1}=-\frac{g}{2}+1$.
\\ Let now be $M>1$ and suppose that
\[ Y(A,P)v_{0}|_{P=P_{+}} = A_{-s_{M}} v_{0} = A \]
where $A_{-s_{M}}=\left( a_{-n_{1}+\frac{g}{2}}....a_{-n_{M}+\frac{g}{2}}\right)$.
\\ Consider now
\[ :\nabla^{k}a(P)Y(A,P): = \sum_{n}:\nabla^{k}a(P)Y(A,P):_{n} f_{M+1}^{n}(P) \]
We have to show that 
\[ :\nabla^{k}a(P)Y(A,P):_{n}v_{0} = 0 \mbox{ for } n> -s_{M+1} \]
in order to get
\[ :\nabla^{k}a(P)Y(A,P):v_{0}|_{P=P_{+}} = :\nabla^{k}a(P)Y(A,P):_{-s_{M+1}} v_{0}. \]
We have by the definition of the normal ordered product:
\[ :\nabla^{k}a(P)Y(Av_{0},P): = \sum_{n}\left( \underbrace{ \sum_{j<\frac{g}{2},m} a^{(k)}_{j}A_{m}l^{jm}_{n} }_{(I)}+ \underbrace{ \sum_{j\geq \frac{g}{2},m} A_{m}a^{(k)}_{j}l^{jm}_{n} }_{(II)} \right) f_{M+1}^{n}(P) \]
where $l^{jm}_{n}= l^{jm}_{n,(M)}$ (see Prop. \ref{prp:nop}).
\\ We get for the sum $(II)$:
\begin{lem} For all $n$ we have 
\[ \sum_{j\geq \frac{g}{2},m} A_{m} a^{(k)}_{j}v_{0}l^{jm}_{n} = 0\]
\end{lem}
This is a consequence of eq. (\ref{eq:zero}) and the condition $a_{\frac{g}{2}+x}v_{0} = 0$ for $x=0,1,2,..$. $\square$
\begin{lem} If $A_{m}v_{0}=0$ for $m>-s_{M}$ then for $n>-s_{M+1}$ we have
\[ \sum_{j<\frac{g}{2},m} a^{(k)}_{j}A_{m}v_{0}l^{jm}_{n} = 0 \]
\end{lem}
Proof. We have to consider $n=-s_{M+1}+x$ where $x=1,2,3,..$ and $j= \frac{g}{2}-y$ where $y=1,2,3,...$. Then we get because of equation (\ref{eq:ljmn}):
\[ l^{\frac{g}{2}-y, m}_{-s_{M+1}+x}\neq 0 \Rightarrow -s_{M+1}-g+x+y \leq m\leq -s_{M+1} +x+y \]
Because of $s_{M+1}=\frac{(-1-2M)g}{2}+M+1 = \frac{(1-2M)g}{2}-g+M+1 = s_{M}-g+1$ 
we have:
\[ l^{\frac{g}{2}-y, m}_{-s_{M+1}+x}\neq 0 \Rightarrow -s_{M}+x+y-1 \leq m\leq -s_{M}+g+x+y-1. \]
It follows that for $x,y=1,2,3,...$ we have $m>-s_{M}$ and because of $A_{m}v_{0}=0$ for $m>-s_{M}$ we get for the first sum: $\sum_{j<\frac{g}{2}} a_{j}^{(k)}A_{m}v_{0}l^{j, m}_{n}=0$. $\square$
\\ Because of $A_{m}v_{0}=0$ for $m>-s_{M}$ we get
\[ \sum_{j<\frac{g}{2},m} a^{(k)}_{j}A_{m}v_{0}l^{jm}_{-s_{M+1}} = (l^{\frac{g}{2}-1, -s_{M}}_{-s_{M+1}}) a_{-k+\frac{g}{2}-1}A_{-s_{M}}v_{0} + ... \hspace{0.5cm} \square \] 
%
%
\section{Translation and Locality Property} 
\subsection{Translation Property}
\begin{prop}[\cite{schsh}] For the field $a(P):= \mathcal{Y}(a_{\frac{g}{2}-1}v_{0},P)=\sum_{n}a_{n}\omega^{n}(P)$ we have
\[ \nabla a(P) = \sum_{n} [T,a_{n}] \omega^{n}(P). \]
\end{prop}
We have furthermore
\[ \nabla^{k+1} a(P)  = [T,\nabla^{k}a(P)] \]
Because of the definition of the normal ordered product we can deduce:
\[ \nabla :\nabla^{k}a(P):\nabla^{h}a(P):: = :\nabla^{k+1}a(P):\nabla^{h}a(P):: + :\nabla^{k}a(P):\nabla^{h+1}a(P): : \]
Altogether we have
\[ \nabla \mathcal{Y}(A,P) = [T,\mathcal{Y}(A,P)]. \]
\subsection{Wick Product and Locality}
We write the field $a(P)=\sum_{n}a_{n}\omega^{n}(P)$ as a sum of two field as follows:
\[ a(P) = a(P)_{+} + a(P)_{-}\]
where $a(P)_{+}=\sum_{n<\frac{g}{2}}a_{n}\omega^{n}(P)$, and $a(P)_{-}=\sum_{n\geq \frac{g}{2}}a_{n}\omega^{n}(P)$.
\\ According to the definition of the normal ordered product (\ref{dfn:nop}) 
for the field $a(P)$ of weight $1$ we get 
\begin{eqnarray*}
 :\nabla^{k}a(P)\nabla^{h}a(P): & = & \nabla^{k}a(P)_{+} \nabla^{h}a(P) + \nabla^{h}a(P) \nabla^{k}a(P)_{-} \\
& = & \nabla^{k}a(P)_{+}\nabla^{h}a(P)_{+} + \nabla^{k}a(P)_{+}\nabla^{h}a(P)_{-} \\
& + & \nabla^{h}a(P)_{+} \nabla^{k}a(P)_{-} + \nabla^{h}a(P)_{-} \nabla^{k}a(P)_{-}  
\end{eqnarray*}
Thus the inductively defined normal ordered product $:\nabla^{n_{1}}a(P)_{+}...\nabla^{n_{M}}a(P)_{+}:$ is a sum of $2^{M}$ summands of the form
\[ \nabla^{n_{1}}a(P)_{+}\nabla^{n_{2}}a(P)_{+}...\nabla^{n_{j-1}}a(P)_{+} \nabla^{n_{M}}a(P)_{-}...\nabla^{n_{j}}a(P)_{-} \]
where $1\leq j\leq M$, that means that after some $j$ the fields $\nabla^{n_{k}}a(P)_{+}$ (for $k>j$) are arranged in the opposite order (see also \cite{kac} section 3.3 for the usual normal ordered product).
\\ We have the relation: $[a_{n},a_{m}]=\gamma_{nm}K$ (see eq. (\ref{eq:gamma})). Therefore we get formally (due to eqs. (\ref{eq:gamma})(\ref{eq:zero})):
\begin{equation} [\nabla^{k}a(P)_{-},\nabla^{h}a(Q)_{-}] = 0 \end{equation}
and
\begin{equation} [ [ \nabla^{k}a(P)_{-}, \nabla^{h}a(Q)],\nabla^{f}a(Q) ] = 0 \mbox{ for all } k,h,f \end{equation}
The last equation follows from the fact that $K$ is the central element.
\\ We can now formulate a Wick product like formula (compare with equation (3.3.3) in \cite{kac}):
\[ \left( :\nabla^{n_{1}} a(P)...\nabla^{n_{M}} a(P): \right)\left(:\nabla^{m_{1}}a(Q)...\nabla^{m_{N}}a(Q):\right) = \]
\[ =  \sum_{s=0}^{\min(M,N)}\sum_{i_{1}<...<i_{s}}^{j_{1}\neq ...\neq j_{s}} [\nabla^{i_{1}}a(P)_{-}, \nabla^{j_{1}}a(Q)]...[\nabla^{i_{s}}a(P)_{-}, \nabla^{j_{s}}a(Q)] \cdot \]
\[ \cdot :\nabla^{k}a(P) ... \nabla^{k}a(P) \nabla^{k}a(Q) ... \nabla^{k}a(Q):_{(i_{1}...i_{s},j_{1}...j_{s})} \]
where the subscript $(i_{1}...i_{s},j_{1}...j_{s})$ means that the fields $\nabla^{i_{1}}a(z)...\nabla^{i_{s}}a(z)$,  and $\nabla^{j_{1}}a(z)...\nabla^{j_{s}}a(z)$ are removed. The normal ordered product $:a(z)b(w):$ for two fields is defined in the obvious way.
\\ If we now consider the "contractions" $[\nabla^{i_{k}}a(P)_{-}, \nabla^{j_{k}}a(Q)]$ as expansions of certain kernels that look locally like
\[ \frac{dzdw}{(z-w)^{h}} +... \mbox{ for some }h\in\Z_{>0} \]
then we have a sort of locality property for our state-field correspondence.
%
%
%
\section{A proposed Global Vertex Algebra}
Inspired by the above considerations we can define a global vertex algebra:
\begin{dfn}
Let be $X$ a Riemann surface, $P_{+}$ and $P_{-}$ two distinguished points on $X$ in general position.
\\ A \textbf{global vertex algebra} is a collection of data:
\\  - A  $\Z_{\geq 0}$ graded vector space $V$, 
\\ - (vacuum vector) a vector $v_{0}\in V$, 
\\ - (translation operator) a linear map $T\in End{V}$,
\\ - (vertex operator) a linear map
\begin{eqnarray*}
V & \rightarrow & End(V)[[(f_{\lambda}^{n}(P))_{\lambda, n\in\Z}]] \\
A & \mapsto & \mathcal{Y}(A,P) = \sum_{n} A_{n}f_{\lambda}^{n}(P).
\end{eqnarray*}
These data are subject to the following axioms:
\begin{enumerate}
\item (vacuum) $\mathcal{Y}(v_{0},P) = id_{V}$. 
\\ For any $A\in V$ of conformal dimension $\lambda$ we have 
\[ \mathcal{Y}(A,P)v_{0} = \sum_{n<-s_{\lambda}}f_{\lambda}^{n}(P). \]
and
\[ \mathcal{Y}(A,P)v_{0}|_{P=P_{0}} = A \]
in other words, $A_{n}v_{0} = 0$, $n> -s_{\lambda}$, and $A_{-s_{\lambda}}v_{0}= A$.
\item (translation) For any $A\in V$, $Tv_{0} = 0$, and
\[ \nabla \mathcal{Y}(A,P) = [T,\mathcal{Y}(A,P)]. \]
\item (locality) 
\[ [\mathcal{Y}(A,P),\mathcal{Y}(B,Q)] (E(P,Q))^{N} \mbox{ holomorphic (but multivalued) on } X\times X.\]
\end{enumerate}
\end{dfn}
\textbf{Remark 7.2} \\
For $g=0$ this definition coincides with the definition of a vertex algebra in the following sense:
\\ The $f_{\lambda}^{n}(P)$ are replaced by the monomials $z^{-n-\lambda}$. 
\\ For the translation axiom $\nabla$ is replaced by $\partial$.
\\ In the locality axiom the Schottky-Klein prime form is replaced by $(z-w)$.
\\ \textbf{Remark 7.3} \\
The definition of a chiral algebra due to Beilinson and Drinfeld \cite{bd} (see also \cite{g} and especially \cite{fre}) is a D-module $\mathfrak{A}$ equipped with a map
\[ \mu : \mathfrak{A}\sqtimes\mathfrak{A}\rightarrow \Delta_{!}\mathfrak{A} = \frac{K\sqtimes \mathfrak{A}(\infty \Delta)}{K\sqtimes \mathfrak{A}}\]
that satisfies certain properties (e.g. skew symmetry and Jacobi identity).
From \cite{fre} (Lemma 18.3.6) we know that the canonical bundle $K$, considered as a D-module (via the Lie derivation), is a chiral algebra, equipped with the map 
\[ \mu_{K}: K\sqtimes K\stackrel{\cong}{\rightarrow} K_{X\times X}\rightarrow \Delta_{!} K. \]
In our situation we have the spaces $\mathcal{F}^{\lambda}$ that become also D-modules via the Lie derivative (see section 2.1).
\\ The map $\mu$ can be regarded as the sheaf theoretic version of the operator product expansion. In our case however the operator product expansion seem to be connected with kernels:
\[ \frac{K^{\lambda_{1}}\sqtimes K^{\lambda_{2}}(\lambda_{1}+\lambda_{2}+2\Delta)} {K^{\lambda_{1}}\sqtimes K^{\lambda_{2}}}. \]
\\ \textbf{Remark 7.4} \\
 It might be interesting to consider this generalization of vertex algebras from a more abstract point of view. The fields in the vertex algebra are formal Laurentseries. $\C((z))$ is a graded algebra by $z^{n}\cdot z^{m} = z^{n+m}$. 
\\ Denote by $U[[\mathfrak{t}_{\lambda}^{n}]]$ a space of infinitely many variables over an associative algebra $U$. Let be a product defined by 
\[ \mathfrak{t}_{\lambda}^{n}\cdot \mathfrak{t}_{\mu}^{m} = \sum_{k=n+m-k_{0}}^{n+m+k_{1}} \eta_{nm}^{k} \mathfrak{t}_{\lambda+\mu}^{k}\mbox{ where } \eta_{nm}^{k}\in U. \]
That means the usual $t^{n}\cdot t^{m}=t^{n+m}$ is replaced by a \textbf{quasigraded} product. An analogous notion for locality might be obtained as well. 
\section{Affine Krichever-Novikov Algebras}
It is tempting to extend the above ideas to Affine Algebras of Krichever-Novikov type.
\\ Here we give only a short description how it could be defined:
\\ Let $\la{g}$ be a complex finite-dimensional reductive Lie algebra. Then
\[ \hat{\la{g}} = \la{g}\otimes \mathcal{A} \oplus \C K\]
is called the affine Lie algebra of KN-type. We have the relations:
\[ [a_{n},b_{m}] = \sum_{k}\alpha_{nm}^{k}[a,b]_{k} + (a|b) \gamma_{nm}K, \hspace{0.3cm} [K,\hat{\la{g}}]= 0, \]
for the definition of $\alpha_{nm}^{k}$ see proposition \ref{res:alpha}, and for $\gamma_{nm}$ see eq. (\ref{eq:gamma}).
\begin{dfn} The delta distribution of weight $\lambda$ is defined by
\[ \Delta_{\lambda}(P,Q) = \sum_{n} f_{\lambda}(P) f^{1-\lambda}(Q). \]
\end{dfn}
These delta distributions on Riemann surfaces can be considered as the formal sum of the expansion of the corresponding Szeg\"{o} kernels.
\\ Set especially $\Delta(P,Q):= \Delta_{0}(P,Q)$.
\begin{prop}[see\cite{schsh}]
\[ d_{P}\Delta(P,Q) = \sum_{nm} \gamma_{n,m} \omega^{n}(P)\omega^{m}(Q). \]
\end{prop}
$d_{P}\Delta(P,Q)$ can be considered as the formal expansion of a differential which looks locally
\[ \frac{dzdw}{(z-w)^{2}}+ ... \]
\begin{cor} We have
\[ [a(P),b(Q)] = [a,b](P)\Delta(P,Q) + d_{P}\Delta(P,Q) \]
and therefore 
\[ [a(P),b(Q)](E(P,Q))^{2} \mbox{ holomorphic on } X\times X. \]
\end{cor}
Proof.
\begin{eqnarray*}
[a(P),b(Q)] & = & \sum_{nm}[a_{n},b_{m}] \omega^{n}(P)\omega^{m}(Q) \\
& = & \sum_{nm} \left(\sum_{k}\alpha_{nm}^{k}[a,b]_{k} + (a|b) \gamma_{nm}K\right)\omega^{n}(P)\omega^{m}(Q) \\ 
& \stackrel{(*)}{=} & \sum_{n} \sum_{k}[a,b]_{k}\omega^{k}(P)A_{m}(P)\omega^{m}(Q) + d_{P}\Delta(P,Q) \\
& = & \sum_{k}[a,b]_{k}\omega^{k}(P) \sum_{m}A_{m}(P)\omega^{m}(Q) + d_{P}\Delta(P,Q)\\
& = & [a,b](P)\Delta(P,Q) + d_{P}\Delta(P,Q)
\end{eqnarray*} 
In $(*)$ we used the symmetry $\alpha_{nm}^{k}=\alpha_{mn}^{k}$, the corollary \ref{res:alpha} and the above proposition. $\square$ \vspace{0.5cm}
\\  \textbf{Acknowledgements}
\\ I want to thank Peter Newstead and Slava Nikulin and Martin Schottenloher for fruitful discussions. I also want to thank the university of Liverpool for hospitality. 
%
\addtocontents{toc}{Bibliography}
%

\end{document}